\documentclass{amsart}
\usepackage{amscd,amssymb,relsize, hyperref, verbatim}
\usepackage[arrow,matrix,graph,frame,poly,arc,tips]{xy}

\newtheorem{thm}{Theorem}[section]
\newtheorem*{theo}{Theorem}
\newtheorem{cor}[thm]{Corollary}

\newtheorem{prop}[thm]{Proposition}

\theoremstyle{definition}
\newtheorem{definition}[thm]{Definition}

\newtheorem{remark}[thm]{Remark}

\numberwithin{equation}{section}

\newcommand{\set}[1]{\left\{#1\right\}}

\renewcommand{\b}[1]{\mathbf{#1}}

\newcommand{\cC}{\mathcal{C}}

\newcommand{\cG}{\mathcal{G}}
\newcommand{\cL}{\mathcal{L}}

\newcommand{\R}{\mathbb{R}}
\newcommand{\C}{\mathbb{C}}
\newcommand{\Q}{\mathbb{Q}}

\newcommand{\Z}{\mathbb{Z}}
\renewcommand{\k}{\Bbbk}

\newcommand{\PS}{{P\!\varSigma}}

\DeclareMathOperator{\id}{id}

\DeclareMathOperator{\Aut}{Aut}

\DeclareMathOperator{\GL}{GL}

\DeclareMathOperator{\Tor}{Tor}
\DeclareMathOperator{\geom dim}{geom\ dim}

\DeclareMathOperator{\tc}{{\sf TC}}
\DeclareMathOperator{\zcl}{\sf zcl}
\DeclareMathOperator{\cat}{cat}
\DeclareMathOperator{\secat}{secat}
\DeclareMathOperator{\cl}{\sf cl}
\DeclareMathOperator{\cd}{cd}

\DeclareMathOperator{\IA}{IA}
\DeclareMathOperator{\Inn}{Inn}

\begin{document}

\title{Topological complexity of basis-conjugating automorphism groups}

\author[D.~C. Cohen]{Daniel C. Cohen$^\dag$}
\address{Department of Mathematics, Louisiana State University,
Baton Rouge, LA 70803}
\email{\href{mailto:cohen@math.lsu.edu}{cohen@math.lsu.edu}}
\urladdr{\href{http://www.math.lsu.edu/~cohen/}
{http://www.math.lsu.edu/\~{}cohen}}
\thanks{{$^\dag$}Partially supported 
by National Security Agency grant H98230-05-1-0055}

\author[G. Pruidze]{Goderdzi Pruidze}
\address{Department of Mathematics, Louisiana State University,
Baton Rouge, LA 70803}
\email{\href{mailto:gio@math.lsu.edu}{gio@math.lsu.edu}}
\urladdr{\href{http://www.math.lsu.edu/~gio/}
{http://www.math.lsu.edu/\~{}gio}}



\begin{abstract}
We compute the topological complexity of Eilenberg-Mac\,Lane spaces associated to the group of automorphisms of a finitely generated free group which act by conjugation on a given basis, and to certain subgroups.
\end{abstract}

\date{April 11, 2008}

\maketitle

\section{Introduction} \label{sec:Intro}
Given a mechanical system, a motion planning algorithm is a function which assigns to any pair of states of the system, an initial state and a desired state, a continuous motion of the system starting at the  initial state and ending at the desired state.  Interest in such algorithms arises in robotics, see Latombe \cite{La} as a general reference.  In a sequence of recent papers 
\cite{fa1}-\cite{fa4}, Farber develops a topological approach to the problem of motion planning, introducing a numerical invariant which gives a measure of the ``navigational complexity'' of the system.

Let $X$ be a path-connected topological space, the space of all possible configurations of a mechanical system.  In topological terms, the motion planning problem consists of finding an algorithm which takes pairs of configurations, i.e., points $(x_0,x_1)\in X\times X$, and produces a continuous path $\gamma\colon [0,1] \to X$ from the initial configuration $x_0=\gamma(0)$ to the terminal configuration $x_1=\gamma(1)$.
Let $PX$ be the space of all continuous paths in $X$, equipped with the compact-open topology.  The map $\pi\colon PX \to X\times X$ defined by sending a path to its endpoints, $\pi\colon \gamma \mapsto(\gamma(0),\gamma(1))$, is a fibration.  The motion planning problem then asks for a section of this fibration, a map 
$s\colon X\times X \to PX$ satisfying $\pi \circ s = \id_{X \times X}$.  It would be desirable for a motion planning algorithm to depend continuously on the input.  However, one can show that there exists a globally continuous motion planning algorithm $s\colon X \times X \to PX$ if and only if $X$ is contractible, see \cite[Thm. 1]{fa1}.  One is thus led to study the discontinuities of such algorithms.

Define the \emph{topological complexity} of $X$ to be the Schwarz genus, or sectional category, of the path-space fibration, 
$\tc(X):=\secat(\pi\colon PX \to X \times X)$.  That is, $\tc(X)$ is the smallest number $k$ for which there is an open cover $X\times X = U_1 \cup \dots \cup U_k$ such that the map $\pi$ admits a continuous section $s_j\colon U_j \to PX$ over each $U_j$, $\pi \circ s_j = \id_{U_j}$. One can show that $\tc(X)$ is an invariant of the homotopy type of $X$, see \cite[Thm. 3]{fa1}.

Let $X$ be an aspherical space, a space whose higher homotopy groups vanish, $\pi_i(X)=0$ for $i\ge 2$.  In \cite[\S31]{fa4}, Farber poses the problem of computing the topological complexity of such a space in terms of algebraic properties of the fundamental group $G=\pi_1(X)$.  In other words, given a discrete group $G$, define the topological complexity of $G$ to be $\tc(G):=\tc(K(G,1))$, the topological complexity of an Eilenberg-Mac{\,}Lane space of type $K(G,1)$, and express $\tc(G)$ in terms of invariants such as the cohomological or geometric dimension of $G$ if possible.  

A number of results in the literature may be interpreted in the context of this problem.  For a right-angled Artin group $G$, the topological complexity of an associated $K(G,1)$-complex was computed in \cite{cp}.  For the Artin pure braid group $G=P_n$, and the group $P_{n,m}=\ker(P_n \to P_m)$, the kernel of the homomorphism which forgets the last $n-k$ strands of a pure braid, the configuration spaces $F(\C,n)$ and $F(\C_m,n)$ of $n$ ordered points in $\C$, and $\C_m=\C\setminus\{m\ \text{points}\}$ respectively are associated Eilenberg-Mac\,Lane spaces.  In \cite{fy} and \cite{fgy}, Farber, Grant, and Yuzvinsky determine the topological complexity of these configuration spaces.  
All of these results may be expressed in terms of the cohomological dimension, $\cd(G)$, of the underlying group $G$.  For instance, one has $\tc(P_n) = \tc(F(\C,n)) = 2n-2=2\cd(P_n)$.

The pure braid group $P_n$ and the group $P_{n,k}$ may be realized as subgroups of 
$\Aut(F_n)$, the automorphism group of the finitely generated free group $F_n=\langle x_1,\dots,x_n\rangle$.  The purpose of this note is to determine the topological complexity of 
several other subgroups of $\Aut(F_n)$.

Let $G=\PS_n$ be the ``group of loops'', the group of motions of a collection of $n\ge 2$ unknotted, unlinked circles in $3$-space, where each (oriented) circle returns to its original position.  This group may be realized as the basis-conjugating automorphism group, or pure symmetric automorphism group, 
of $F_n$, consisting of all automorphisms which, for the fixed basis $\{x_1,\dots,x_n\}$ for $F_n$, send each generator to a conjugate of itself.  
A presentation for $\PS_n$ was found by McCool \cite{mc}.  In particular, this group is generated by automorphisms $\alpha_{i,j} \in \Aut(F_n)$, $1\le i \neq j
\le n$, defined by $\alpha_{i,j}(x_i)=x_j^{} x_i^{} x_j^{-1}$ and $\alpha_{i,j}(x_k)=x_k$ for $k \neq i$.
Also of interest is the ``upper triangular McCool group,'' 
the subgroup $\PS_n^+$ of $\PS_n$ generated by $\alpha_{i,j}$ for $i<j$.  The main results of this note may be summarized as follows.

\begin{theo} \label{thm:main}
For the basis-conjugating automorphism group $\PS_n$ and the upper triangular McCool group $\PS_n^+$, the topological complexity is given by
\[
\tc(\PS_n) = 2n-1\quad\text{and}\quad \tc(\PS_n^+)=2n-2.
\]
\end{theo}

Let $X$ be an Eilenberg-Mac\,Lane complex of type $K(G,1)$ for either $G=\PS_n$ or $G=\PS_n^+$.  Since the topological complexity of $X$, $\tc(X)=\tc(G)$, is the Schwarz genus of the path-space fibration, it admits several useful bounds.  
For instance, one has
\[
\tc(X)=\secat(\pi\colon PX \to X \times X) \le \cat(X \times X) \le 2\cat(X)-1 \le 2\dim(X)+1,
\]
where $\cat(X)$ denotes the Lusternik-Schnirelmann category of $X$, 
see Schwarz \cite{Schwarz} and James \cite{Ja} as classical references.  One also has 
a cohomological lower bound 
\[
\tc(X) \ge 1+ \cl(\ker(\pi^*\colon 
H^*(X\times X;\Q)\to H^*(PX;\Q))),
\]
where $\cl(A)$ denotes the cup length of a graded ring $A$, the largest integer $q$ for which there are homogeneous elements $a_1,\dots,a_q$ of positive degree in $A$ such that $a_1\cdots a_q \neq 0$.  Using the K\"unneth formula, the fact that $PX \simeq X$, and the equality 
$H^*(X;\Q)=H^*(G;\Q)$, the kernel of $\pi^*\colon 
H^*(X\times X;\Q)\to H^*(PX;\Q)$ may be identified with the kernel $Z=Z(H^*(G;\Q))$ of the cup-product map $H^*(G;\Q) \otimes H^*(G;\Q) \xrightarrow{\cup} H^*(G;\Q)$, see \cite[Thm. 7]{fa1}.  The cup length of the ideal $Z$ of zero-divisors is referred to as the \emph{zero-divisor cup length} of $H^*(G;\Q)$, and is denoted by $\zcl(H^*(G;\Q))=\cl(Z)$.  In this notation, the cohomological lower bound reads
\[
\tc(G) \ge 1+ \zcl(H^*(G;\Q)).
\]

This note is organized as follows.  After a discussion of basis-conjugating automorphism groups in Section \ref{sec:bcg}, including the determination of their geometric dimensions, we use the (known) structure of the cohomology rings of these groups to compute the zero-divisor cup lengths of these rings in Section \ref{sec:coh}.  These results are used in Section \ref{sec:tc} to find the topological complexity of these groups.  
We conclude with some remarks concerning formality in Section~\ref{sec:fm}.

\section{ Basis-conjugating automorphism groups} \label{sec:bcg}

Let $N$ be a compact set contained in the interior of a manifold $M$.  Generalizing the familiar interpretation of a braid as the motion of $N=\set{n\ \text{distinct points}}$ in $M=\R^2$, Dahm \cite{d} defines a motion of $N$ in $M$ as a path $h_t$ in $\mathcal{H}_c(M)$, the space of homeomorphisms of $M$ with compact support, satisfying $h_0=\id_M$ and $h_1(N)=N$.  With an appropriate notion of equivalence, the set of equivalence classes of motions of $N$ in $M$ is a group, and, furthermore, there is a homomorphism from this group to the automorphism group of the fundamental group $\pi_1(M\setminus N)$.

In \cite{gold}, Goldsmith gives an exposition of Dahm's (unpublished) work, with particular attention paid to the case where $N=\cL_n$ is a collection of $n$ unknotted, unlinked circles in $M=\R^3$.  Let $\cG_n$ denote the corresponding motion group.  Goldsmith shows that $\cG_n$ is generated by three types of motions, flipping a single circle, interchanging two (adjacent) circles, and pulling one circle through another, and that the Dahm homomorphism $\phi\colon\cG_n \to \Aut(\pi_1(\R^3\setminus \cL_n))$ is an embedding.

Fix a basepoint $e\in\R^3$ disjoint from $\cL_n=C_1 \cup\dots\cup C_n$, and for each $i$, let $x_i$ be (the homotopy class of) a loop based at $e$ linking $C_i$ once.  This identifies $\pi_1(\R^3\setminus\cL_n,e)=F_n$ with the free group generated by $x_1,\dots,x_n$.  With this identification, the generators of the motion group $\cG_n\hookrightarrow\Aut(F_n)$ correspond to automorphisms $\rho_i$ (flip $C_i$), $\tau_i$ (switch $C_i$ and $C_{i+1}$), and $\alpha_{i,j}$ (pull $C_i$ through $C_j$) defined by
\begin{equation*} \label{eq:autos}
\rho_{i}(x_k) =
\begin{cases}
x_{k}^{-1} & \text{if $k= i$,}\\
x_k & \text{if $k \neq i$,}
\end{cases}
\qquad
\tau_{i}(x_k) =
\begin{cases}
x_{k+1} & \text{if $k = i$,}\\
x_{k-1} & \text{if $k = i+1$,}\\
x_k & \text{if $k \neq i, i+1$,}
\end{cases}
\end{equation*}
and
\begin{equation} \label{eq:alpha}
\alpha_{i,j}(x_k) =
\begin{cases}
x_{j}x_{k}x_{j}^{-1} & \text{if $k =i$,}\\
x_k & \text{if $k \neq i$.}
\end{cases}
\end{equation}

Let $\varphi: \Aut(F_{n}) \rightarrow \Aut(F_{n}/ [F_{n}, F_{n} ]) \cong \GL(n,\Z)$ denote the epimorphism induced by the abelianization homomorphism $F_n \to F_n/[F_n,F_n]\cong \Z^n$. There is a corresponding short exact sequence 
$1 \rightarrow \IA_{n}  \rightarrow \Aut(F_{n}) \xrightarrow{\varphi} \GL(n,\Z) \rightarrow 1$, 
where $\IA_{n}=\ker \varphi$ is the well known group of automorphisms of $F_n$ which induce the identity on $H_1(F_n;\Z)$.  
In \cite{bl}, Brownstein and Lee considered the following commutative diagram 
\begin{equation*}  \label{eq:diagram}
 \begin{CD}
  1 @>>> \ker(\varphi \circ \phi) @>>>  \cG_{n}           @>\varphi \circ \phi>> \Z /2 \wr \Sigma_{n} @>>>1 \\
   @.                 @VVV                           @V\phi VV                                                 @VVV                      @.\\
  1 @>>> \IA_{n}                         @>>> \Aut(F_{n})     @>>> \GL(n, \Z)                                            @>>> 1
 \end{CD}
\end{equation*}
where the vertical maps are embeddings, and showed that the image of $\cG_{n}$ under $\varphi\circ\phi$ is the wreath product $ \Z /2 \wr \Sigma_{n}$, the reflection group of type $\mathrm{D}_{n}$. The kernel of $\varphi \circ \phi$ corresponds to the group $\cC_n$ of ``pure motions'' of $\cL_{n}$, motions which bring each oriented circle back to its original position. The isomorphic image of $\ker (\varphi \circ \phi)$ in $\Aut (F_{n})$, i.e., the intersection $\IA_{n}  \cap \phi (\cG_{n})$, is the basis-conjugating automorphism group of the free group.

\begin{definition} \label{def:bcg}
The \emph{basis-conjugating automorphism group} of the free group $F_n$ is the subgroup of $\Aut(F_{n})$ generated by the elements $\alpha_{i,j}$ from \eqref{eq:alpha} with $ 1\le i,j \le n$, and $i\neq j$. Following \cite{jmm}, we denote this group by $\PS_{n}$. 
\end{definition} 

In \cite{mc}, McCool showed that $\PS_n$ admits a presentation with the aforementioned generators, and defining relations
\begin{equation} \label{eq:BCrels}
\left\{\begin{matrix}
 [\alpha_{i,j},\alpha_{k,l}] \hfill &
  \text{for $i,j,k,l$ distinct} \hfill  \\
[\alpha_{i,j},\alpha_{k,j}]  \hfill\ &
\text{for $i,j,k$ distinct} \hfill \\
[\alpha_{i,j},\alpha_{i,k}\alpha_{j,k}] \hfill &
\text{for $i,j,k$ distinct} \hfill
\end{matrix}\right\},
\end{equation}
where $[\alpha, \beta]=\alpha^{}\beta^{}\alpha^{-1}\beta^{-1}$ denotes the commutator.

An ``upper triangular'' version of the basis-conjugating automorphism group has been an object of study in a number of recent works, see \cite{BM,ccp,cpvw}.

\begin{definition} \label{def:utriang}
The \emph{upper triangular McCool group} $\PS_{n}^{+}$ is the subgroup of $\PS_{n}$ generated by the elements $\alpha_{i,j}$ with $i<j$, subject to the relevant relations \eqref{eq:BCrels}.
\end{definition}

The upper triangular McCool group $\PS_n^+$ shares a number of features with the Artin pure braid group $P_n$.  For instance, both groups may be realized as iterated semidirect products of free groups:
\[
P_n = F_{n-1} \rtimes_{\rho_{n-1}} \cdots \rtimes_{\rho_2} \rtimes F_1
\quad\text{and}\quad
\PS_n^+ = F_{n-1} \rtimes_{\mu_{n-1}} \cdots \rtimes_{\mu_2} \rtimes F_1.
\]
For the pure braid group, the action of the free group $F_k$ on $F_m$ with $1\le k<m\le n-1$ is given by the restriction of the Artin representation $\rho_m\colon P_m \to \Aut(F_m)$, see for instance \cite{Birman}.
For the upper triangular McCool group, the action of $F_k=\langle \alpha_{n-k,j}\mid n-k+1\le j\le n\rangle$ on $F_m=\langle \alpha_{n-m,j}\mid n-m+1\le j\le n\rangle$, that is, the homomorphism 
$\mu_m\colon \rtimes_{j=1}^{m-1} F_j \to \Aut(F_m)$, was determined in \cite{cpvw} (with different notation).  Using the relations \eqref{eq:BCrels}, one can check that
\[
\mu_{m}(\alpha_{j,p})(\alpha_{i,q})=
\alpha_{j,p}^{-1}\alpha_{i,q}^{}\alpha_{j,p}^{}=
\begin{cases}
\alpha_{i,p}^{} \alpha_{i,q}^{} \alpha_{i,p}^{-1}&\text{if $q=j$,}\\
\alpha_{i,q}^{}&\text{otherwise,}
\end{cases}
\]
where $i=n-m$, $j=n-k$, $1\le i < j< p \le n$, and $i+1\le q \le n$.

Consideration of centers provides another similarity between these groups.  For a group $G$, let $Z(G)$ denote the center of $G$, and let $\overline{G}=G/Z(G)$.  It is well known that the center of the pure braid group is infinite cyclic, and that $P_n \cong \overline{P}_{n} \times Z(P_n)=\overline{P}_n \times \Z$.  The analogous result holds for the upper triangular McCool group.

\begin{prop} \label{prop:plus product}
The center of the upper triangular McCool group $\PS_n^+$ is infinite cyclic, the quotient 
$\overline{\PS}_n^+ =F_{n-1} \rtimes_{\mu_{n-1}} \cdots \rtimes_{\mu_{3}} F_2$ is an iterated semidirect product of free groups, and $\PS_n^+ \cong\overline{\PS}_n^+ \times Z(\PS_n^+)=\overline{\PS}_n^+ \times \Z$. 
\end{prop}
\begin{proof}
Consider the element $c=\alpha_{1,n} \alpha_{2,n}\cdots \alpha_{n-1,n}$ of the group $\PS^{+}_{n}$.  Using \eqref{eq:BCrels}, it is readily checked that $c$ commutes with all the generators of $\PS^{+}_{n}$, so $c\in Z(\PS^{+}_{n})$. Furthermore, it is clear that $c \in \Aut(F_n)$ has infinite order. Consequently, the infinite cyclic subgroup $C=\langle c\rangle$ is contained in the center $Z(\PS_n^+)$.

Since $\alpha_{n-1,n} = (\alpha_{1,n}\alpha_{2,n}\cdots\alpha_{n-2,n})^{-1}\cdot c$, the group $\PS_n^+$ admits a presentation with generators $c$ and $\alpha_{i,j}$, $1\le i<j\le n$ and $(i,j)\neq (n-1,n)$, and relations $[c,\alpha_{i,j}]$ for all $i<j$, and the relations \eqref{eq:BCrels} (not involving $\alpha_{n-1,n}$).  Thus, $\PS_n^+ \cong C \times (\PS_n^+/C)$.  Since the free group $F_1$ in the iterated semidirect product decomposition $\PS_n^+=\rtimes_{j=1}^{n-1} F_j$ is generated by $\alpha_{n-1,n}$, it is clear from the above discussion that $\PS_n^+/C=F_{n-1} \rtimes_{\mu_{n-1}} \cdots \rtimes_{\mu_{3}} F_2$.  An easy inductive argument reveals that the center of this quotient is trivial.  It follows that $C=Z(\PS_n^+)$, which completes the proof.
\end{proof}

Despite the aforementioned similarities, the groups $P_n$ and $\PS_n^+$ are not isomorphic, see Bardakov and Mikhailov \cite{BM}.

\begin{definition} \label{def:cohdim}
Let $G$ be a group. The \emph{cohomological dimension} of $G$, denoted by $\cd(G)$, is the smallest integer $n$ such that $H^q(G;M)=0$ for any $G$-module $M$ and all $q>n$. 
The \emph{geometric dimension} of the group $G$, denoted by $\geom dim(G)$, is the smallest dimension of an Eilenberg-Mac\,Lane complex of type $K(G,1)$.
\end{definition}

\begin{prop}\label{prop:dimp}
Let $\PS_{n}$ be the basis-conjugating automorphism group. Then
\[
\geom dim (\PS_{n})=\cd (\PS_{n})= n-1. 
\]
\end{prop}

\begin{proof}
As shown by Collins \cite{coll}, for each $n$, the cohomological dimension of $\PS_n$ is as asserted, $\cd(\PS_n)=n-1$.  By a classical result of Eilenberg and Ganea \cite{eg}, for groups of cohomological dimension at least $3$, the geometric dimension is equal to the cohomological dimension.  
Thus, the assertion holds for $\PS_n$ with $n \ge 3$.

Since $\PS_2=F_2$ is the free group generated by $\alpha_{2,1}$ and 
$\alpha_{1,2}$, the case $n=2$ is immediate.

It remains to consider the case $n=3$.  The group $\PS_3$ is generated by six elements $\alpha_{i,j}$, $1\le i\neq j \le 3$.  Let 
$\beta_1=\alpha_{2,1}\alpha_{3,1}$,
$\beta_2=\alpha_{1,2}\alpha_{3,2}$, and 
$\beta_3=\alpha_{1,3}\alpha_{2,3}$, and observe that 
these elements generate the inner automorphism group $\Inn(F_3)$ of $F_3$, which is isomorphic to $F_3$.
As noted in \cite{bl}, 
the group $\PS_3=\Inn(F_3) \rtimes F$ is a semidirect product, where $F=\langle
\alpha_{1,2},\alpha_{2,1},\alpha_{3,1}\rangle$ is also a free group on $3$ generators.  Thus, $\PS_3 \cong F_3 \rtimes F_3$ is a semidirect product of two finitely generated free groups.  

In \cite[\S1.3]{cs1}, Cohen and Suciu give an explicit construction of a $K(G,1)$-complex $X_G$ for an arbitrary iterated semidirect product of finitely generated free groups $G$.  If $G=\rtimes_{i=1}^\ell F_{d_i}$, the complex $X_G$ is $\ell$-dimensional.  In particular, for the group $G=\PS_3$, this construction yields a $2$-dimensional $K(G,1)$-complex.  We therefore have 
$\geom dim (\PS_{3})=\cd (\PS_3) = 2$.
\end{proof}

A similar result holds for the upper triangular McCool groups. 
\begin{prop} \label{prop:dimg}
Let $\PS_{n}^{+}$ be the upper triangular McCool group, and $\overline{\PS}_n^+=\PS_n^+/Z(\PS_n^+)$. Then
\[
\geom dim (\PS_{n}^{+})=\cd (\PS_{n}^{+})= n-1\ \text{and}\ 
\geom dim (\overline{\PS}_{n}^{+})=\cd (\overline{\PS}_{n}^{+})= n-2.
\]
\end{prop}
\begin{proof}
Since $\overline{\PS}_n^+=F_{n-1} \rtimes_{\mu_{n-1}}  \cdots \rtimes_{\mu_3} \rtimes F_2$ 
and $\PS_n^+=\overline{\PS}_n^+ \times \Z$ are iterated semidirect products of finitely generated free groups, this follows immediately from the results of \cite{cs1}.
\end{proof}

\section{Structure of the cohomology ring} \label{sec:coh}

As noted in the Introduction, the zero-divisor cup length of the cohomology ring of a group provides a lower bound for the topological complexity.  In this section, we determine this lower bound for the groups $\PS_{n}$, and 
$\PS_{n}^{+}$.

Let $A=\bigoplus_{k=0}^\ell A^k$ be a graded algebra over a field $\k$, and recall that the cup length $\cl(A)$ is the largest integer $q$ for which there are homogeneous elements $a_1,\dots,a_q$ of positive degree in $A$ such that $a_1\cdots a_q \neq 0$.  The tensor product $A \otimes A$ has a natural graded algebra structure, with multiplication $(u_1\otimes v_1)\cdot (u_2\otimes v_2)=
(-1)^{|v_1|\cdot|u_2|} u_1 u_2 \otimes v_1 v_2$.  Let $\mu\ \colon A \otimes A \to A$ denote the multiplication homomorphism, and let $Z=\ker(\mu)$ be the ideal of zero-divisors.  The zero-divisor cup length of $A$, denoted by $\zcl(A)$, is the cup length of this ideal, $\zcl(A)=\cl(Z)$.  Observe that if $a \in A$, then the element $\bar{a}=a\otimes 1-1\otimes a \in Z$ is a zero-divisor.

In \cite{bl}, Brownstein and Lee determined the low-dimensional cohomology $H^{\le 2}(\PS_n;\Z)$ of the basis-conjugating automorphism group, and conjectured the general ring structure in terms of generators and relations.  This conjecture was recently proved by Jensen, McCammond, and Meier \cite[Thm.~6.7]{jmm}.  For our purposes, it suffices to work with coefficients in the field $\k=\Q$ of rational numbers.  So we  suppress coefficients, and denote the rational cohomology of a group $G$ by  $H^*(G)=H^*(G;\Q)$ throughout this section and the next.

\begin{thm}[\cite{jmm}] \label{thm:bt1}
The rational cohomology algebra $H^*(\PS_n)$ is isomorphic to $E/I$, where $E$ is the exterior algebra over $\Q$ generated by degree one elements $a_{i,j}$, $1\le i \neq j \le n$, and $I$ is the homogeneous ideal generated by the degree two elements
\[
a_{i,j} a_{j,i}\ (i,j\ \text{distinct})\quad\text{and}\quad a_{k,j} a_{j,i} - a_{k,j}a_{k,i} - a_{i,j} a_{k,i}\ (i,j,k\ \text{distinct}).
\]
\end{thm}

This result may be used to exhibit an explicit basis for $H^q(\PS_n)$ for each $q$, $0\le q \le n-1$, see \cite[\S6]{jmm}.  Call an element of the form $a_{i,j}a_{j,k} \cdots a_{s,t}a_{t,i}$ a cyclic product.  Then, $H^q(\PS_n)$ has a basis consisting of those $q$-fold products $a_{i_1,j_1}a_{i_2,j_2}\cdots a_{i_q,j_q}$ of the one-dimensional generators which do not contain any cyclic products, and have distinct first indices $i_1,\dots,i_q$.  It follows that the Poincar\'e polynomial of $\PS_n$ is $\sum_{q \ge 0} \dim H^q(\PS_n)\cdot t^q = (1+nt)^{n-1}$.  In particular, $H^i(\PS_n)=0$ for $i\ge n$, and the cup length of $H^*(\PS_n)$ is $n-1$.

These results may be used to determine the zero-divisor cup length of the ring $H^*(\PS_n)$.

\begin{thm} \label{thm:s1}
Let $\PS_{n}$ be the basis-conjugating automorphism group. Then the zero-divisor cup length of the rational cohomology algebra of $\PS_n$ is 
\[
\zcl( H^*(\PS _{n}))  = 2n-2.
\]
\end{thm}
\begin{proof}
In general, the zero-divisor cup length of an algebra $A$ cannot exceed the cup length of the tensor product $A\otimes A$, which is twice the cup length of $A$ itself, $\zcl(A) \le \cl(A \otimes A) = 2\cl(A)$.  Since $\cl(H^*(\PS_n))=n-1$ by Theorem \ref{thm:bt1}, it follows that $\zcl(H^*(\PS_n))\le 2n-2$.  

To establish the reverse inequality, we work in the aforementioned basis for $H^*(\PS_n)$, and the corresponding induced basis for the tensor product $H^*(\PS_n)\otimes H^*(\PS_n)$.
Observe that any monomial in the generators of $H^*(\PS_n)$ that contains a cyclic product must vanish, and that any finite expression in $H ^* (\PS _n ) $ can be reduced to an expression in the basis elements after finitely many applications of the relation
\begin{equation} \label{eq:BCcohrel}
  a_{k,j}a_{k,i} = a_{k,j}a_{j,i}+a_{i,j}a_{k,i},
\end{equation} 
step-by-step eliminating repetition in the first index. 

For each $i < n$, consider the elements $\b{x}_{i}=a_{i,i+1}$ and $\b{y}_{i}=a_{i+1,i}$ in $H^{*}(\PS_{n})$, and the corresponding zero divisors  $\bar{\b{x}}_{i} = \b{x}_i \otimes 1 - 1 \otimes \b{x}_i$ and $\bar{\b{y}}_{i} = \b{y}_{i} \otimes 1 - 1 \otimes \b{y}_{i}$ in the tensor product $ H^{*}(\PS_{n}) \otimes H^{*}(\PS_{n})$.
We claim that the product 
\[
M= \prod_{i=1}^{n-1} \bar{\b{x}}_i \cdot \prod_{i=1}^{n-1} \bar{\b{y}}_i =
\bar{\b{x}}_1\bar{\b{x}}_2\cdots \bar{\b{x}}_{n-1}\bar{\b{y}}_1\bar{\b{y}}_2\cdots \bar{\b{y}}_{n-1}
\]
of these $2n-2$ zero divisors is different from zero. To prove this claim, we use the relation \eqref{eq:BCcohrel} to express $M$ in terms of the specified basis of the tensor product, and identify at least one monomial that stays unaffected throughout the reduction process.

If $I$ is a subset of $[n-1]=\{1,2,\dots,n-1\}$, let $|I|$ denote the cardinality of $I$, and let $U_{I}=z_1 \cdots z_{n-1}$ and $V_{I}= \hat{z}_1 \cdots \hat{z}_{n-1}$,  where 
\[
z_{i} =
\begin{cases}
\b{y}_{i}, & \text{if $i \notin I$, } \\
\b{x}_{i}, & \text{if $i \in I$ }
\end{cases}
\text{\ and\ }
\hat{z}_{i} =
\begin{cases}
\b{y}_{i}, & \text{if $i \in I$, } \\
\b{x}_{i}, & \text{if $i \notin I$. }
\end{cases}
\]
Then, using the fact that $\bar{\b{x}}_{i} \bar{\b{y}}_{i} = \b{y}_{i} \otimes \b{x}_{i} - \b{x}_{i} \otimes \b{y}_{i}$, one has
\begin{equation} \label{eq:M}
M=\sum_{I \subseteq [1,n-1] } (-1)^{\left| I\right|} U_{I}\otimes V_{I}. 
\end{equation}
When $I=\emptyset$ is the empty set, 
the summand $U_{\emptyset} \otimes V_{\emptyset}$ in \eqref{eq:M} is 
\[
U_{\emptyset} \otimes V_{\emptyset} = \b{y}_{n-1}\b{y}_{n-2}\cdots \b{y}_{1} \otimes \b{x}_{1}\b{x}_{2} \cdots \b{x}_{n-1} = a_{n,n-1}a_{n-1,n-2} \cdots a_{2,1} \otimes a_{1,2}a_{2,3} \cdots a_{n-1,n}.
\]
This monomial is already a basis element of $H^{n-1} (\PS_{n}) \otimes H^{n-1}(\PS_{n})$. 

We claim that the expression of any other summand $(-1)^{|I|}U_I \otimes V_I$ of 
\eqref{eq:M} in terms of our basis for $H^*(\PS_n)\otimes H^*(\PS_n)$ will avoid 
the specified basis element $U_{\emptyset} \otimes V_{\emptyset}$.  Clearly, if the monomial $U_I$ is already a basis element of $H^*(\PS_n)$, there is nothing to prove.
Otherwise, $U_{I}$ contains a factor $a_{k,j}a_{k,i}$ for at least one $k$ with $ 1 < k < n $, and these are the only generators in the product 
$U_I$ involving index $k$. 
Applying the relation \eqref{eq:BCcohrel} to the product $a_{k,j}a_{k,i}$, we obtain (up to sign)
\[
U_{I} =( a_{k,j}a_{j,i}+a_{i,j}a_{k,i} ) \cdot \{ \text{other factors} \} = a_{k,j} P + a_{k,i} Q, 
\]
where $P$ and $Q$ are monomials in the generators $a_{r,s}$ of $H^*(\PS_n)$ 
with $r\neq k$ and $s\neq k$.
Further application of reductive relation \eqref{eq:BCcohrel} to $P $ and $Q$ will result in no further appearance of $k$ in the indices. 
Consequently, expressing $U_I=a_{k,j} P + a_{k,i} Q$ in the specified basis for $H^*(\PS_n)$ will yield a linear combination of basis elements, each with exactly one factor involving  index $k$.  On the other hand, our fixed monomial 
$U_{\emptyset} = a_{n,n-1}\cdots a_{k+1,k} a_{k,k-1} \cdots a_{2,1} $ contains two factors involving index $k$. Hence, the basis monomial $U_{\emptyset} \otimes V_{\emptyset} $ is different from any other possible basis summand coming from $U_{I}\otimes V_{I} $ with $I \ne \emptyset$, and our claim holds.
\end{proof}

The cohomology of the upper-triangular McCool group $\PS_n^+$ may be analyzed in a similar manner.  The integral cohomology of $\PS_n^+$ was computed by Cohen, Pakianathan, Vershinin, and Wu \cite[Thm.~1.4]{cpvw}.  Their results imply the following.

\begin{thm}[\cite{cpvw}] \label{thm:bt2}
The rational cohomology algebra $H^*(\PS_n^+)$ is isomorphic to $E^+/I^+$, where $E^+$ is the exterior algebra over $\Q$ generated by degree one elements $a_{i,j}$, $1\le i < j \le n$, and $I^+$ is the homogeneous ideal generated by the degree two elements 
\[
a_{i,j} a_{i,k} - a_{i,j} a_{j,k}\ (i < j < k).
\]
\end{thm}

This result may be used to exhibit an explicit basis for $H^q(\PS_n^+)$ for each $q$, $0\le q \le n-1$, compare \cite[\S7]{cpvw}.   The group $H^q(\PS_n^+)$ has a basis consisting of those $q$-fold products $a_{i_1,j_1}a_{i_2,j_2}\cdots a_{i_q,j_q}$ of the one-dimensional generators which satisfy
$1\le i_1<i_2<\cdots <i_q\le n-1$ and $i_p < j_p \le n$ for each $p$.  It follows that 
$\sum_{q \ge 0} \dim H^q(\PS_n^+)\cdot t^q = \prod_{k=1}^{n-1}(1+kt)$.  
In particular, $H^i(\PS_n^+)=0$ for $i\ge n$, and the cup length of $H^*(\PS_n^+)$ is $n-1$.

These results facilitate analysis of the zero-divisor cup length of the ring $H^*(\PS_n^+)$.

\begin{thm} \label{thm:s2}
Let $\PS_{n}^+$ be the upper-triangular McCool group. Then the zero-divisor cup length of the rational cohomology algebra of $\PS_n^+$ satisfies
\[
\zcl(H^{*}(\PS^{+}_{n})) \ge 2n-3.
\]
\end{thm}
\begin{proof}
The basis for $H^*(\PS_n^+)$ specified above induces a basis for $H^*(\PS_n^+)\otimes H^*(\PS_n^+)$.

We check that the product of the $2n-3$ zero-divisor elements 
\begin{equation} \label{eq:N}
\bar{a}_{1,n-1} \bar{a}_{1,n} \bar{a}_{2,n-1} \bar{a}_{2,n}\cdots \bar{a}_{n-2,n-1} \bar{a}_{n-2,n}\cdot (a_{n-1,n} \otimes a_{n-1,n} ) 
\end{equation}
is nonzero, where $\bar{a}_{i,j} = a_{i,j}\otimes 1 - 1 \otimes a_{i,j}$. Note that 
\[
\bar {a}_{i,n-1} \cdot  \bar{a}_{i,n} =     a_{i,n} \otimes a_{i,n-1} -    a_{i,n-1} \otimes a_{i,n}   + a_{i,n-1}a_{i,n} \otimes 1  +1\otimes a_{i,n-1}a_{i,n} 
\]
for any $i\le n-2$. The product \eqref{eq:N} contains summands  of the form
\begin{equation} \label{eq:survive}
a_{1,i_1} a_{2,i_2} \cdots a_{n-2,i_{n-2}}a_{n-1,n} \otimes a_{1,j_{1}} a_{2,j_{2}} \cdots a_{n-2,j_{n-2}} a_{n-1,n},
\end{equation}
where $i_{p}$ and $ j_{p}$ take different values from the set $ \{n-1,n \}$ for each $p$. Such summands represent different basis elements in the tensor product. Observe that these are the only terms which have exactly one factor in the both sides of the tensor product with first subindex $q$, for each $q$, $1 \le q \le n-1$. In light of the relations $a_{i,j} a_{i,k} = a_{i,j} a_{j,k}$ in $H^*(\PS_n^+)$, expressing other summands in \eqref{eq:N} in terms of the specified basis for 
 $H^*(\PS_n^+)\otimes H^*(\PS_n^+)$ cannot yield terms with this feature.  
Thus, the terms given by \eqref{eq:survive} survive, and the product \eqref{eq:N} is nonzero.
\end{proof}

\begin{remark} It follows from the results of the next section that equality holds in Theorem~\ref{thm:s2}, $\zcl(H^*(\PS_n^+))=2n-3$.
\end{remark}

\section{Topological complexity} \label{sec:tc} 
In this section, we recall several necessary properties of topological complexity, and prove the main results of the paper.

Let $X$ be a path-connected topological space.  We are interested in the the case where $X$ is an Eilenberg-Mac\,Lane space of type $K(G,1)$ for $G=\PS_n$ or $G=\PS_n^+$, so assume that $X$ has the homotopy type of a finite CW-complex.  Let $PX$ denote the space of all continuous paths 
$\gamma\colon [0,1] \to X$, equipped with the compact-open topology.  The map $\pi\colon PX \to X \times X$, $\gamma \mapsto (\gamma(0),\gamma(1))$, defined by sending a path to its endpoints is a fibration, with fiber $\Omega{X}$, the based loop space of $X$.  

Recall from the Introduction that the motion planning problem asks for a (continuous) section of this fibration, a map  $s\colon X\times X \to PX$ satisfying $\pi \circ s = \id_{X\times X}$.  As shown by Farber \cite[Thm.~1]{fa1}, in most cases, such a section cannot exist.

\begin{prop}[\cite{fa1}] \label{prop:contr}
The path space fibration $\pi: PX \to X\times X$ admits a section if and only if $X$ is contractible.
\end{prop}

\begin{definition}
The \emph{topological complexity} of $X$, $\tc(X)$,  is the smallest positive integer $k$ for which $X\times X=U_1\cup\dots\cup U_k$, where $U_j$ is open and there exists a continuous section $s_j\colon U_i\to PX$, $\pi \circ s_j=\id_{U_i}$, for each $j$, $1\le j \le k$.  In other words, the topological complexity of $X$ is
the Schwarz genus (or sectional category) of the path space fibration $\pi\colon PX \to X\times X$.
\end{definition}

The topological complexity of $X$ is a homotopy-type invariant, see \cite[Thm.~3]{fa1}.  If $G$ is a discrete group, define $\tc(G)$, the topological complexity of $G$, to be that of an Eilenberg-Mac\,Lane space of type $K(G,1)$.  In \cite[\S31]{fa4}, Farber poses the problem of determining the topological complexity of $G$ in terms of other invariants of $G$, such as $\cd(G)$, the cohomological dimension.  In this section, we solve this problem for the basis-conjugating automorphism groups 
$\PS_n$ and $\PS_n^+$.

We will require several properties of topological complexity.  We briefly record these, and refer to the survey \cite{fa4} for further details.

First, if $X$ is a finite-dimensional cell complex, then $\tc(X) \le 2 \dim(X)+1$, see \cite[\S3]{fa4}.  
Consequently, if $G$ is a group of finite geometric dimension, then
\begin{equation} \label{eq:upper}
\tc(G) \le 2 \geom dim(G)+1.
\end{equation}
Second, as noted in the Introduction, a lower bound for the topological complexity of a group $G$ is provided by the zero-divisor cup length of the cohomology ring $H^*(G)=H^*(G;\Q)$:
\begin{equation} \label{eq:lower}
\tc(G) \ge 1+\zcl(H^*(G)),
\end{equation}
see \cite[\S15]{fa4}.  Finally, if $X $ and $Y$ are path-connected paracompact locally contractible topological spaces (in particular, CW-complexes), then $\tc (X \times Y)  \le \tc(X) + \tc(Y) -1$, 
see \cite[\S12]{fa4}.  
Consequently, if $G_1$ and $G_2$ are groups (of finite geometric dimension), then
\begin{equation} \label{eq:tcprod}
\tc (G_1 \times G_2)  \le \tc(G_1) + \tc(G_2) -1.
\end{equation}

With these facts at hand, we now prove our main theorems. 

\begin{thm} \label{thm:m1}
The topological complexity of the basis-conjugating automorphism group $\PS_n$~is
\[
\tc(\PS_n)=2n-1.
\]
\end{thm}
\begin{proof}
By Theorem \ref{thm:s1}, the zero-divisor cup length of $H^*(\PS_n)$ is given by $\zcl(H^*(\PS_n))=2n-2$.  So the lower bound \eqref{eq:lower} yields $\tc(\PS_n) \ge 2n-1$.  For the reverse inequality, recall from Proposition \ref{prop:dimp} that $\geom dim(\PS_n)=\cd(\PS_n)=n-1$. Consequently, the upper bound \eqref{eq:upper} yields $\tc(\PS_n) \le 2n-1$, completing the proof. 
\end{proof}

\begin{thm} \label{thm:m2}
The topological complexity of the upper triangular McCool group $\PS_n^+$ is
\[
\tc (\PS^{+}_{n}) =2n-2.
\]
\end{thm}

\begin{proof}
By Theorem \ref{thm:s2}, the zero-divisor cup length of $H^*(\PS_n^+)$ satisfies 
$\zcl(H^*(\PS_n^+))\ge 2n-3$.  So the lower bound \eqref{eq:lower} yields $\tc(\PS_n^+) \ge 2n-2$.  

For the reverse inequality, recall from Proposition \ref{prop:plus product} that $\PS_n^+ \cong \overline{\PS}_n^+ \times \Z$.  Since the circle $S^1$ is a $K(\Z,1)$-space, and $\tc(\Z)=\tc(S^1)=2$ (see, for instance, \cite[\S5]{fa1}), the product inequality \eqref{eq:tcprod} yields
\[
\tc(\PS_n^+) \le \tc(\overline{\PS}_n^+) + \tc(\Z) - 1 = \tc(\overline{\PS}_n^+) + 1.
\]
By Proposition \ref{prop:dimg}, we have $\geom dim(\overline{\PS}_n^+)=\cd(\overline{\PS}_n^+)=n-2$.  Consequently, the upper bound \eqref{eq:upper} yields $\tc(\overline{\PS}_n^+)\le 2n-3$. 
Thus, $\tc(\PS_n^+)\le 2n-2$, completing the proof.
\end{proof}

\begin{cor} The zero-divisor cup length of the rational cohomology algebra $H^*(\PS_n^+)$ is 
$\zcl(H^*(\PS_n^+))=2n-3$.
\end{cor}

\section{Formality} \label{sec:fm}

If $X$ is an Eilenberg-Mac\,Lane space of type $K(G,1)$, where $G=\PS_n$ or $G=\PS_n^+$, the results of the previous section imply that the topological complexity of $X$ is given by the cohomological lower bound, 
\[
\tc(X) = 1+\zcl(H^*(X;\Q)).
\]
This equality holds for a number of spaces of interest in topology, including certain configuration spaces, complements of certain complex hyperplane arrangements, and Eilenberg-Mac\,Lane spaces corresponding to right-angled Artin groups, see \cite{cp, fgy,fy,yuz}.  Since all of these spaces are formal in the sense of Sullivan \cite{sul}, it is natural to speculate that such an equality holds for an arbitrary formal space $X$, conjecturally, $\tc(X)=1+\zcl(H^*(X;R))$ for appropriate coefficients $R$.  This is explicitly conjectured for the complement of an arbitrary hyperplane arrangement by Yuzvinsky in \cite{yuz}.  Related problems are studied in \cite{FGKV,lm}.  
In this section, we show that the upper triangular McCool group $\PS_n^+$ provides evidence in favor of such a conjecture.

\begin{thm} \label{thm:form}
Let $X$ be an Eilenberg-Mac\,Lane space of type $K(G,1)$, where $G=\PS_n^+$ is the upper triangular McCool group.  Then $X$ is a formal space.
\end{thm}

In order to prove this theorem, we will need some definitions and facts concerning formality and related notions.

Let $X$ be a space with the homotopy type of a connected, finite-type CW-complex.  Loosely speaking, $X$ is \emph{formal} if the rational homotopy type of $X$ is determined by the rational cohomology ring $H^{*}(X; \Q)$. Examples of formal spaces include spheres, simply-connected Eilenberg-Mac\,Lane spaces, and  those mentioned above. 

Let $G$ be a finitely presented group.  Following Quillen \cite{Q}, call $G$ 
\emph{$1$-formal} if the Malcev Lie algebra of $G$ is 
quadratic, see \cite{ps2} for details.  As shown by Sullivan \cite{sul} and Morgan \cite{Mo}, the fundamental group $G=\pi_1(X)$ of a formal space $X$ is a $1$-formal group.  There are, however, non-formal spaces with $1$-formal fundamental groups, see \cite{K,Mo}.

Papadima and Suciu \cite[Prop.~2.1]{ps} provide a sufficient condition for the formality of a CW-complex.  Recall that a connected, graded algebra $A$ over a field $\k$ is said to be a \emph{Koszul algebra} if $\Tor_{p,q}^A(\k,\k)=0$ for all $p \neq q$, where $p$ is the homological degree of the $\Tor$ groups, and $q$ is the internal degree coming from the grading of $A$.  A necessary condition is that $A$ be a quadratic algebra, the quotient of a free algebra on generators in degree $1$ by an ideal generated in degree $2$.

\begin{prop}[\cite{ps}]
\label{prop:ps}
Let $X$ be a connected, finite-type CW-complex. If $H^{*}(X; \Q)$ is a Koszul algebra and $G=\pi_{1}(X)$ is a $1$-formal group, then $X$ is a formal space.
\end{prop}

Berceanu and Papadima \cite[Rem.~5.5]{bp} have recently shown that the upper triangular McCool group $\PS_n^+$ is $1$-formal.  Thus, to prove Theorem \ref{thm:form}, it suffices to show that the rational cohomology algebra $H^*(\PS_n^+;\Q)$ is Koszul.  For this, we 
will use a result of Jambu and Papadima \cite[Prop.~6.3]{jp}.

Let $A=\bigoplus_{k\ge 0} A^k$ be a connected, graded $\k$-algebra, and denote the augmentation ideal of $A$ by $A^+=\bigoplus_{k\ge 1}A^k$. Call a subalgebra $B$ of $A$ normal if $AB^+ = B^+ A$.  If $B \subset A$ is normal, there is a canonical projection $\pi \colon A \to F$, where $F=A/AB^+$.

\begin{prop}[\cite{jp}] 
\label{prop:jambupapad}
Let $B\subset A$ be a normal subalgebra such that $A$ is free as a right $B$-module, and assume that the $\k$-algebras $A$, $B$ and $F= A / AB^+$ are quadratic. If $B$ and $F$ are Koszul algebras, then $A$ is a Koszul algebra.
\end{prop}

We apply this result to the rational cohomology algebra $H^*(\PS^{+}_{n}; \Q)$.

\begin{prop} \label{prop:formality}
The rational cohomology algebra $H^*(\PS^{+}_{n}; \Q)$ of the upper triangular McCool group is a Koszul algebra.
\end{prop}

\begin{proof}  Write $A_n = H^*(\PS_n^+;\Q)$.  

The proof consists of an inductive application of Proposition \ref{prop:jambupapad}. 
Since $\PS_2^+ \cong \Z$, the base case $A_2$ is trivial.

Inductively assume that $A_{n-1}$ is Koszul.  For $k<n$, observe that $A_k$ is isomorphic to the subalgebra $\tilde{A}_k$ of $A_n$ generated by the elements $a_{i,j}$ with $n-k<i<j\le n$.  Thus, we may assume that the subalgebra $\tilde{A}_{n-1}$ of $A_n$ is Koszul.  Since the algebras under consideration are graded commutative, 
$\tilde{A}_{n-1}$ is a normal subalgebra of $A_{n}$. Furthermore, $A_n$ is free as a right $\tilde{A}_{n-1}$-module. Namely, 
\[
   A_{n} = 1 \cdot \tilde{A}_{n-1} \oplus a_{1,2}\cdot \tilde{A}_{n-1} \oplus \cdots \oplus a_{1,n}\cdot
\tilde{A}_{n-1}.
\]
This follows from the fact that in any monomial of the algebra $A_n$, the factor $a_{1,i}$ with minimal $i$ always survives, since $a_{1,i}a_{1,j}=a_{1,i}a_{i,j}$ 
in $A_n$ for any $1<i<j$, see Theorem \ref{thm:bt2}.

Analyzing again the relations in $A_n$, we observe that the algebra $A_n / A_n  \tilde{A}^{+}_{n-1}$ is a graded algebra generated by the elements $a_{1,i}$, $2\le i \le n$, where all the terms in degree $2$ and higher die. Consequently, the algebra $A_n / A_n  \tilde{A}^{+}_{n-1}$ is quadratic, and moreover, Koszul. Thus, all the algebras under consideration are quadratic.  Therefore, the conditions of Proposition \ref{prop:jambupapad} are satisfied, and the result follows immediately.
\end{proof}

Since the upper triangular McCool group $\PS_n^+$ is $1$-formal \cite{bp} and $H^*(\PS_n^+;\Q)$ is Koszul, Proposition \ref{prop:ps} implies that an Eilenberg-Mac\,Lane space of type $K(\PS_n^+,1)$ is formal, proving Theorem \ref{thm:form}.  Such a space $X$ provides an example of a non-simply-connected formal space with $\tc(X)=1+\zcl(H^*(X;\Q))$.

\begin{remark}
Berceanu and Papadima \cite[Thm.~5.4]{bp} also showed that the basis-conjugating automorphism group $\PS_n$ is $1$-formal. We do not know if the cohomology algebra  $H^*(\PS_n;\Q)$ is Koszul.
\end{remark}

\newcommand{\arxiv}[1]{{\texttt{\href{http://arxiv.org/abs/#1}{{arxiv:#1}}}}}
\newcommand{\MRh}[1]{\href{http://www.ams.org/mathscinet-getitem?mr=#1}{MR#1}}

\bibliographystyle{amsplain}

\end{document}